\documentclass[
    ,final            
  ]
  {aipproc}

\layoutstyle{8x11single}

  \usepackage{latexsym}
  \usepackage{amsmath}   
  \usepackage{amsfonts}
  \usepackage{amssymb}
  \usepackage{graphicx}
  \usepackage{amsthm}

  \newcommand{\field}[1]{\mathbb{#1}}
  
  \newcommand{\R}{\field{R}}

  \newcommand{\vect}[1]{\ensuremath{\mbox{\textbf{\textit{#1}}}}}
  \newcommand{\svect}[1]{\ensuremath{\mbox{\textbf{\textit{\small #1}}}}}

  \newcommand{\be}{\begin{equation}}
  \newcommand{\ee}{\end{equation}} 
  \newcommand{\bv}[1]{\mathbf{#1}}
  \newcommand{\blue}[1]{{#1}}

  \newcommand{\bfr}{\begin{frame}}
  \newcommand{\efr}{\end{frame}}

\begin{document}

\title{Angles between subspaces computed in Clifford Algebra}


\author{Eckhard Hitzer}{
  address={Department of Applied Physics, University of Fukui, 910-8507 Japan}
}

\begin{abstract}
We first review the definition of the angle between subspaces and how it is computed using matrix algebra. Then we introduce the Grassmann and Clifford algebra description of subspaces. The geometric product of two subspaces yields the full relative angular information in an explicit manner. We explain and interpret the result of the geometric product of subspaces gaining thus full practical access to the relative orientation information. 
\end{abstract}

\maketitle



\begin{figure}[b]
  \resizebox{0.15\textwidth}{!}{\includegraphics{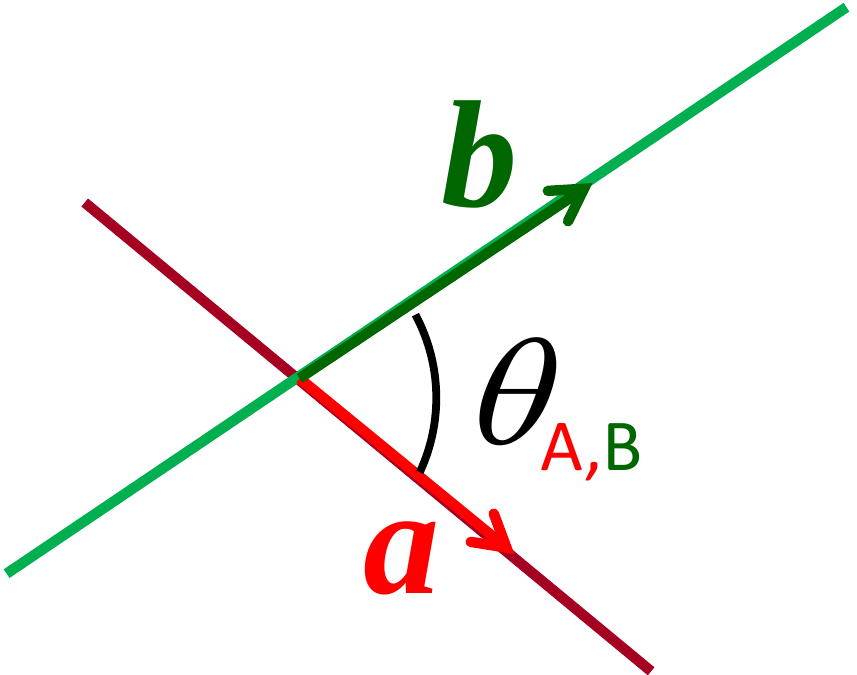}}
  \resizebox{0.4\textwidth}{!}{\includegraphics{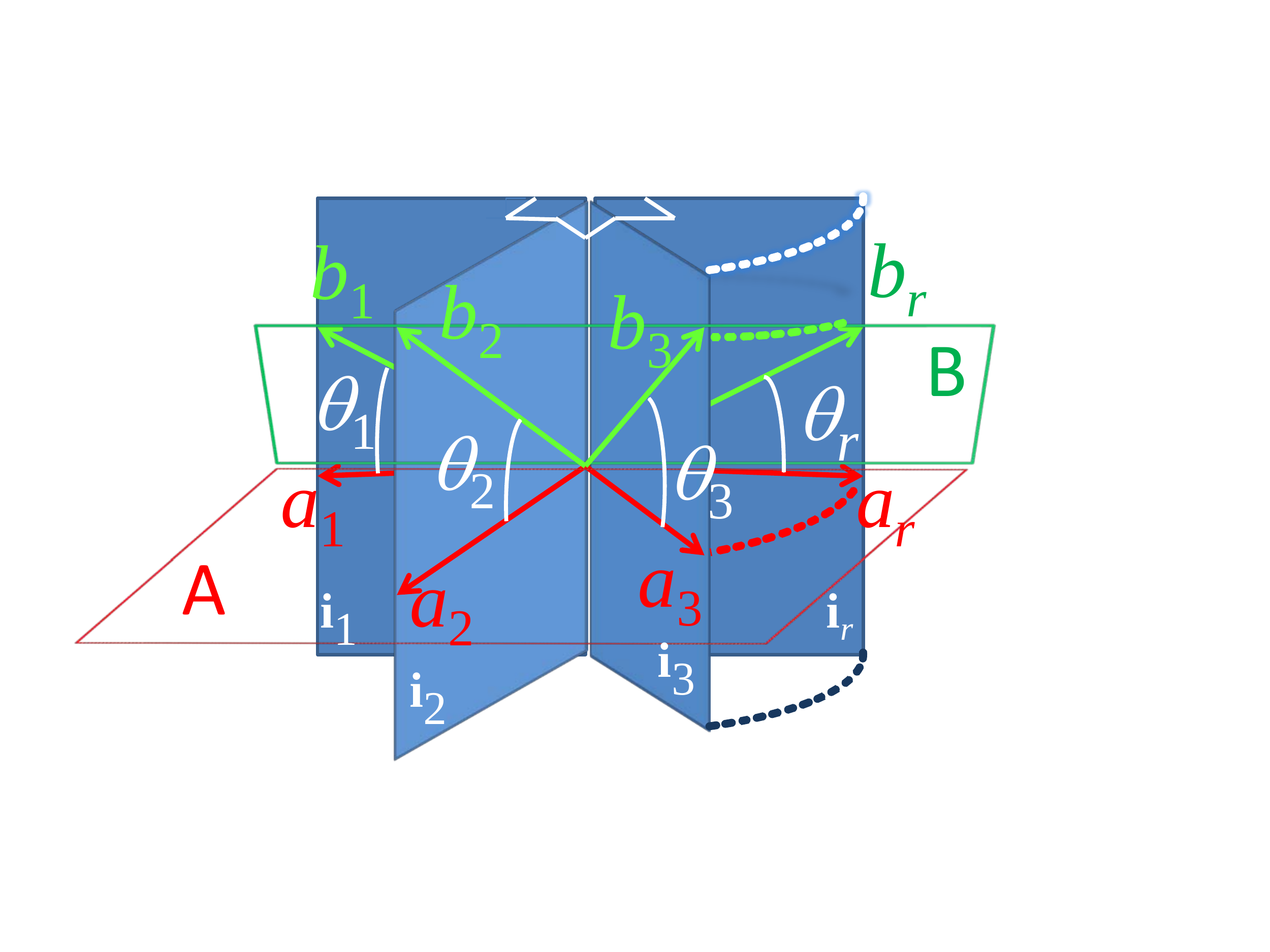}}
  \resizebox{0.4\textwidth}{!}{\includegraphics{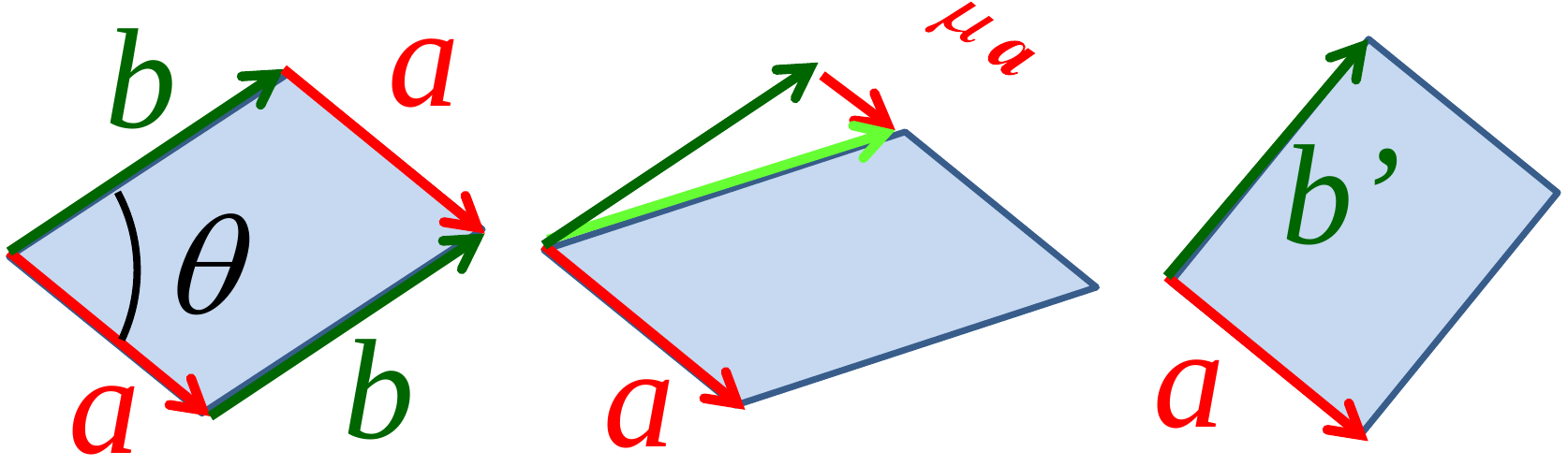}}
  \caption{Left: Angle $\theta_{\mathsf{A},\mathsf{B}}$ between two lines $\mathsf{A}$, $\mathsf{B}$, spanned by unit vectors $\vect{a}$, $\vect{b}$, respectively. 
Center: Angular relationship of two subspaces $\mathsf{A}$, $\mathsf{B}$, spanned by two sets of vectors $\{\vect{a}_1,\ldots,\vect{a}_r\}$, and $\{\vect{b}_1,\ldots,\vect{b}_r\}$, respectively.
Right: Bivectors $\vect{a}\wedge \vect{b}$ as oriented area elements can be freely reshaped in the same plane (e.g. by $\vect{b} \rightarrow \vect{b}+\mu \vect{a}, \mu \in \R$) without changing their value (area and orientation). Orthogonal reshaping into the form of an oriented rectangle is shown. 
  \label{fg:2subspaces}}
\end{figure}

To begin with let us look (see Fig. \ref{fg:2subspaces}, left) at two lines $\mathsf{A}, \mathsf{B}$ in a vector space $\R^n$, which are spanned by two (unit) vectors $\vect{a},\vect{b} \in \R^n, \vect{a}\cdot\vect{a}= \vect{b}\cdot \vect{b} = 1$: 
\be 
  \mathsf{A} = \text{span}[\vect{a}], \quad
  \mathsf{B} = \text{span}[\vect{b}]. 
\ee 
The angle $0\leq \theta_{\mathsf{A},\mathsf{B}}\leq \pi/2$ between lines $\mathsf{A}$ and $\mathsf{B}$ is simply given by
$ 
  \cos \theta_{\mathsf{A},\mathsf{B}} = \vect{a}\cdot\vect{b}.
$


Next let us examine the case of two $r$-dimensional $(r\leq n)$ subspaces $\mathsf{A}, \mathsf{B}$ of an $n$-dimensional Euclidean vector space $\R^n$. The situation is depicted in Fig. \ref{fg:2subspaces}, center. 
Each subspace $\mathsf{A}$, $\mathsf{B}$ is spanned by a set of $r$ linearly independent vectors
\begin{align} 
  \mathsf{A} = \text{span}[\vect{a}_1,\ldots,\vect{a}_r] \subset \R^n, 
   \qquad
  \mathsf{B} = \text{span}[\vect{b}_1,\ldots,\vect{b}_r] \subset \R^n. 
\end{align}
Using Fig. \ref{fg:2subspaces} (center) we introduce the following notation for principal vectors. The angular relationship between the subspaces $\mathsf{A}$, $\mathsf{B}$ is characterized by a set of $r$ principal angles $\theta_k, 1 \leq k \leq r $, as indicated in Fig. \ref{fg:2subspaces} (center). A principal angle is the angle between two principal vectors $\vect{a}_k\in \mathsf{A}$ and  $\vect{b}_k\in \mathsf{B}$. The spanning sets of vectors $\{\vect{a}_1,\ldots,\vect{a}_r\}$, and $\{\vect{b}_1,\ldots,\vect{b}_r\}$ can be chosen such that pairs of vectors $\vect{a}_k, \vect{b}_k$ either 
agree $\vect{a}_k = \vect{b}_k$, $\theta_k=0$, or enclose a finite angle 
$0 < \theta_{k}\leq \pi/2$. 

In addition the pairs of vectors $\{\vect{a}_k, \vect{b}_k\}, 1 \leq k \leq r $ span mutually orthogonal lines (for $\theta_k=0$) and (principal) planes $\bv{i}_k$ (for $0 < \theta_{k}\leq \pi/2$). These mutually orthogonal planes $\bv{i}_k$ are indicated in Fig. \ref{fg:2subspaces} (center). Therefore if $\vect{a}_k \nparallel \vect{b}_k$ and $\vect{a}_l \nparallel \vect{b}_l$ for $1 \leq k\neq l \leq r$, then plane $\bv{i}_k \perp \text{plane } \bv{i}_l$. The cosines of the \textit{principal} angles $\theta_k$ may therefore be $\cos \theta_k = 1$ (for $\vect{a}_k = \vect{b}_k$), or $\cos \theta_k = 0$ (for $\vect{a}_k \perp \vect{b}_k$), or any value $0<\cos \theta_k < 1$. The total angle between the two subspaces $\mathsf{A}$, $\mathsf{B}$ is defined as the product
\be 
  \cos \theta_{\mathsf{A},\mathsf{B}}
  = \cos \theta_1 \cos \theta_2 \ldots \cos \theta_r.
\ee 
In this definition $\cos \theta_{\mathsf{A},\mathsf{B}}$ will automatically be zero if any pair of principal vectors $\{\vect{a}_k, \vect{b}_k\}, 1 \leq k \leq r $ is perpendicular. Then the two subspaces are said to be perpendicular $\mathsf{A} \perp \mathsf{B}$, a familiar notion from three dimensions, where two perpendicular planes $\mathsf{A}$, $\mathsf{B}$ share a common line spanned by $\vect{a}_1 = \vect{b}_1$, and have two mutually orthogonal principal vectors $\vect{a}_2 \perp \vect{b}_2$, which are both in turn orthogonal to the common line vector $\vect{a}_1 $. It is further possible to choose the indexes of the vector pairs $\{\vect{a}_k, \vect{b}_k\}, 1 \leq k \leq r $ such that the principal angles $\theta_k$ appear ordered by magnitude
$
  \theta_1 \geq \theta_2 \geq \ldots \geq \theta_r. 
$


The conventional (matrix algebra) method of computing the angle $\theta_{\mathsf{A},\mathsf{B}}$ between two $r$-dimensional subspaces $\mathsf{A},\mathsf{B} \subset \R^n$ spanned by two sets of vectors $\{\vect{a}'_1,\vect{a}'_2, \ldots \vect{a}'_r \}$ and $\{\vect{b}'_1,\vect{b}'_2, \ldots \vect{b}'_r \}$ is to first arrange these vectors as column vectors into two $n\times r$ matrices
$
  M_\mathsf{A} = [\vect{a}'_1,\ldots,\vect{a}'_r], 
  M_\mathsf{B} = [\vect{b}'_1,\ldots,\vect{b}'_r]. 
$
Then standard matrix algebra methods of QR decomposition and singular value decomposition are applied to obtain 
$r$ pairs of singular unit vectors $\vect{a}_k, \vect{b}_k$ and
$r$ singular values $\sigma_k = \cos \theta_k = \vect{a}_k \cdot \vect{b}_k$. 
This approach uses $O(r^3)$ floating point multiplications.

\section{Even more subtle ways}

Per Ake Wedin in his 1983 contribution \cite{PAW:ABS} to a conference on Matrix Pencils entitled \textit{On Angles between Subspaces of a Finite Dimensional Inner Product Space} first carefully treats the above mentioned matrix algebra approach to computing the angle $\theta_{\mathsf{A},\mathsf{B}}$ in great detail and clarity. Towards the end of his paper he dedicates less than one page to mentioning an alternative method starting out with the words: \textit{But there are even more subtle ways to define angle functions.} 

There he essentially reviews how $r$-dimensional subspaces $\mathsf{A},\mathsf{B} \subset \R^n$ can be represented by $r$-vectors (blades) in Grassmann algebra $A,B \in \Lambda(\R^n)$:
\begin{align} 
  \mathsf{A} = \{\vect{x} \in \R^n | x\wedge A = 0\}, 
  \mathsf{B} =  \{\vect{x} \in \R^n | x\wedge B = 0\}. 
\end{align}
The angle $\theta_{\mathsf{A},\mathsf{B}}$ between the two subspaces $\mathsf{A},\mathsf{B} \in \R^n$ can then be computed in a single step (see also \cite{LD:GAfCS}, p. 68)
\be  
  \cos \theta_{\mathsf{A},\mathsf{B}}
  = \frac{A\cdot \widetilde{B}}{|A| |B|}
  = \cos \theta_1 \cos \theta_2 \ldots \cos \theta_r, 
\label{eq:cosAB}
\ee
where the inner product is canonically defined on the Grassmann algebra $\Lambda(\R^n)$ corresponding to the geometry of $\R^n$. The tilde operation is the reverse operation representing a dimension dependent sign change
$\widetilde{B} = (-1)^{\frac{r(r-1)}{2}}B$,
and $|A|$ represents the norm of blade $A$, i.e. 
$|A|^2 = A \cdot \widetilde{A}$, and similarly $|B|^2 = B \cdot \widetilde{B}$. Wedin refers to earlier works \cite{LA:ABS,QKL:EGoES}. 

Yet equipping a Grassmann algebra $\Lambda(\R^n)$ with a canonical inner product comes close \cite{HL:IAaGR} to introducing Clifford's geometric algebra $Cl_n = Cl(\R^n)$. We will exploit the full geometric product of two subspaces, not only the inner product part \eqref{eq:cosAB}, in order to gain all relative orientation information. And we will treat the case of two subspaces of different dimensions. To the authors knowledge these details have not been published in the literature so far. 

\section{Information in geometric product of two subspace blades}

Clifford (geometric) algebra is based on the geometric product of vectors $\vect{a},\vect{b} \in \R^{p,q}, p+q=n$ 
\begin{equation}
  \vect{a}\vect{b} = \vect{a}\cdot\vect{b} + \vect{a}\wedge\vect{b},
\end{equation}
and the associative algebra $Cl_{p,q}$ thus generated with $\R$ and $\R^{p,q}$ as subspaces of $Cl_{p,q}$. $\vect{a}\cdot\vect{b}$ is the symmetric inner product of vectors and $ \vect{a}\wedge\vect{b}$ is Grassmann's outer product of vectors representing the oriented parallelogram area spanned by $\vect{a},\vect{b}$, compare Fig. \ref{fg:2subspaces} (right).

A \blue{blade} $D_r=\vect{b}_1\wedge\vect{b}_2\wedge\ldots\wedge\vect{b}_r, \vect{b}_l\in\R^{n}, 1\leq l \leq r \leq n$ describes an $r$-dimensional vector \textit{subspace} 
  \be
    \mathsf{D}=\{ \vect{x}\in \R^{p,q} | \vect{x}\wedge D =0 \}.
  \ee 
The magnitude of the blade $D_r \in Cl_n$ is nothing but the volume of the $r$-dimensional parallelepiped spanned by the vectors 
$\{\vect{b}_1,\vect{b}_2,\ldots,\vect{b}_r\}$. 

The geometric product of two $r$-blades $A,B$ contains therefore at most the following grades
\be 
  AB =
  \langle AB\rangle_0 + \langle AB\rangle_2 + \ldots 
  +\langle AB\rangle_{2\text{min}(r,\lfloor n/2\rfloor)},
\ee 
where the limit $\lfloor n/2\rfloor$ (entire part of $n/2$) is due to the dimension limit of $\R^n$.

Every $r$-blade $A_r$ can therefore be written (use e.g. Gram-Schmidt orthogonalization \cite{HS:CAtoGC}) as a product of the scalar magnitude $|A_r| $ times the geometric product of exactly $r$ mutually orthogonal unit vectors $\{\widehat{\vect{a}}_1, \ldots, \widehat{\vect{a}}_r\}$
\be 
  A_r = |A_r| \widehat{\vect{a}}_1 \widehat{\vect{a}}_2 \ldots \widehat{\vect{a}}_r.
\ee 
Please note well, that this \textit{rewriting} of an $r$-blade in geometric algebra does not influence the overall result on the left side, the $r$-blade $A_r$ is before and after the rewriting the very same element of the geometric algebra $Cl_n$. But for the geometric interpretation of the geometric product $AB$ of two $r$-blades $A,B\in Cl_n$ the orthogonal reshaping is indeed a key step. 

A rotation operator (rotor) is given by
$ \vect{x} \longrightarrow (\vect{a}\vect{b})\vect{x}\vect{a}\vect{b}=R^{-1} \vect{x} R,
R=\vect{a}\vect{b} \propto \cos \theta_{\svect{a},\svect{b}}+\bv{i}_{\svect{a},\svect{b}}\sin \theta_{\svect{a},\svect{b}}$, where the unit bivector $\bv{i}_{\svect{a},\svect{b}}$ represents the plane of rotation, and $\theta_{\svect{a},\svect{b}}$ half the rotation angle.


From the discussion of $r$-dimensional subspaces $\mathsf{A}, \mathsf{B}\subset \R^n$ represented by the blades $A=\vect{a}'_1\wedge \ldots \vect{a}'_r$ and $B=\vect{b}'_1\wedge \ldots \vect{b}'_r$, from the freedom of orthogonally reshaping these blades and factoring out the blade magnitudes $|A|$ and $|B|$, and from the classical results of matrix algebra, we now know that we can rewrite the geometric product $AB$ in mutually orthogonal products of pairs of principal vectors $\vect{a}_k,\vect{b}_k, 1\leq k \leq r$
\begin{align} \label{eq:ABrots}
  A\widetilde{B}
  = \vect{a}_1\vect{a}_2\ldots \vect{a}_r\vect{b}_r\ldots \vect{b}_2\vect{b}_1
  = \vect{a}_1\vect{b}_1\vect{a}_2\vect{b}_2\ldots \vect{a}_r\vect{b}_r.
\end{align}
The geometric product 
\begin{align} 
\vect{a}_r\vect{b}_r
  = |\vect{a}_r||\vect{b}_r|(\cos \theta_{\svect{a}_r,\svect{b}_r} +\bv{i}_{\svect{a}_r,\svect{b}_r}\sin\theta_{\svect{a}_r,\svect{b}_r})
  = |\vect{a}_r||\vect{b}_r| (c_r +\bv{i}_rs_r), 
\end{align}
with $c_r=\cos \theta_{\svect{a}_r,\svect{b}_r}$ and
 $s_r=\sin \theta_{\svect{a}_r,\svect{b}_r}$
in the above expression for  $A\widetilde{B}$ is composed of a scalar and a bivector part. The latter is proportional to the unit bivector $\bv{i}_r$ representing a (principal) plane orthogonal to \textit{all} other principal vectors. $\bv{i}_r$ therefore commutes with all other principal vectors, and hence the whole product $\vect{a}_r\vect{b}_r$ (a rotor) commutes with all other principal vectors. A completely analogous consideration applies to all products of pairs of principal vectors, which proves the second equality in \eqref{eq:ABrots}. 

We thus find that we can always rewrite the product $A\widetilde{B}$ as a product of rotors
\begin{align} 
  A\widetilde{B}
  &= |A||B|(c_1 +\bv{i}_1s_1)(c_2 +\bv{i}_2s_2)\ldots(c_r +\bv{i}_rs_r)
   \nonumber \\
  &= |A||B|(c_1c_2\ldots c_r+ 
\,\,\,s_1c_2\ldots c_r\bv{i}_1+c_1s_2\ldots c_r\bv{i}_2+\ldots+c_1c_2\ldots s_r\bv{i}_r+
  \,\,\,\ldots\,
  +s_1s_2\ldots s_r\bv{i}_1\bv{i}_2\ldots\bv{i}_r)
  \label{eq:ABrrspec}
\end{align}

We realize how the scalar part 
$\langle A\widetilde{B} \rangle_0=|A||B|c_1c_2\ldots c_r$, 
the bivector part 
$\langle A\widetilde{B} \rangle_2=|A||B|(s_1c_2\ldots c_r\bv{i}_1+c_1s_2\ldots c_r\bv{i}_2+\ldots+c_1c_2\ldots s_r\bv{i}_r)$, 
etc., up to the $2r$-vector (or $2\lfloor n/2\rfloor$-vector) part 
$\langle A\widetilde{B} \rangle_{2r}=|A||B|s_1s_2\ldots s_r\bv{i}_1\bv{i}_2\ldots\bv{i}_r$ of the geometric product $A\widetilde{B}$ arise and what information they carry. 

Obviously the scalar part yields the cosine of the angle \eqref{eq:cosAB} between the subspaces represented by the two $r$-vectors $A,B\in Cl_n$
\be 
  \cos \theta_{\mathsf{A}\mathsf{B}} 
  = \frac{\langle A\widetilde{B} \rangle_0}{|A||B|}
  = \frac{A\ast\widetilde{B}}{|A||B|}.
\ee

The bivector part consists of a sum of (principal) plane bivectors, which can in general be uniquely decomposed  into its constituent
sum of 2-blades by the method of Riesz, described also in \cite{HS:CAtoGC},  chapter 3-4, equation (4.11a) and following. 

The magnitude of the $2r$-vector part allows to compute the product of all sines of principal angles
\be 
  s_1s_2\ldots s_r = \pm \frac{|\langle A\widetilde{B} \rangle_{2r}|}{|A||B|}. 
\ee

Next, let us now consider the \textit{full relative angular information} between two general $r$-dimensional subspaces $\mathsf{A},\mathsf{B}$, which we take to partly intersect and to be partly perpendicular. We mean by that, that the dimension of the intersecting subspace be $s\leq r$ ($s$ is therefore the number of principal angles equal zero), and the number of principal angles with value $\pi/2$ be $t\leq r-s$. For simplicity we now work with normed blades (i.e. after dividing with 
$|A||B|$). The geometric product of the the $r$-blades $A,B\in Cl_n$ then takes the form (compare the part in brackets with \eqref{eq:ABrrspec})
\begin{align} 
  A\widetilde{B}
  =& (c_{s+1}c_{s+2}\ldots c_{r-t} \,\,\,+ \,\,\,
  s_{s+1}c_{s+2}\ldots c_{r-t}\bv{i}_{s+1}
        +c_{s+1}s_{s+2}\ldots c_{r-t}\bv{i}_{s+2}+
   \ldots
        +c_{s+1}c_{s+2}\ldots s_{r-t}\bv{i}_{r-t}
  \nonumber \\
  &+
 \ldots+\,\,\,s_{s+1}s_{s+2}\ldots s_{r-t}\bv{i}_{s+1}\bv{i}_{s+2}\ldots\bv{i}_{r-t})
  \,\bv{i}_{r-t+1} \ldots \bv{i}_{r}.
\label{eq:ABrrgen}
\end{align}
We thus see, that apart from the integer dimensions $s$ for parallelity (identical to the dimension of the meet of blade $A$ with blade $B$) and $t$ for perpendicularity, the lowest non-zero grade of dimension $2t$ gives the relevant angular measure
$ 
  \cos\theta_{\mathsf{A}\mathsf{B}}=\cos \theta_{s+1}\cos \theta_{s+2}\ldots \cos \theta_{r-t}.
$
While the maximum grade part gives again the product of the corresponding sinuses
$ 
  \sin \theta_{s+1}\sin \theta_{s+2}\ldots \sin \theta_{r-t}.
$

Dividing the product $A\widetilde{B}$ by its lowest grade part  
$\cos\theta_{\mathsf{A}\mathsf{B}} \bv{i}_{r-t+1} \ldots \bv{i}_{r}$ gives a multivector with maximum grade $2(r-t-s)$, scalar part one, and bivector part 
\be 
t_{s+1}\bv{i}_{s+1}
 +t_{s+2} \bv{i}_{s+2}+\ldots
 +t_{r-t}\bv{i}_{r-t},
\ee
where $t_k = \tan \theta_k $. 
Splitting this bivector (see chp. 3-4 of \cite{HS:CAtoGC}) into its constituent bivector parts further yields the (principal) plane bivectors and the tangent values of the principal angles $\theta_k, s<k\leq r-t$.

We finally consider the case of two subspaces $\mathsf{A}, \mathsf{B}$ of \textit{different dimension}. They shall be represented by the normed blades $A=\langle A\rangle_{r+q}$, $r+q, q\geq 0$ and $B=\langle B\rangle_r$. We formally assume the factorized form $A=A_1A_2$, $A_1=\langle A_1\rangle_{q}$, $A2=\langle A_2\rangle_{r}$, with $A_1 \perp A_2$ and  $A_1 \perp B$, 
i.e. $\forall \vect{a} : \vect{a}\rfloor A_1 = \vect{a}A_1$ we have $\vect{a}\rfloor B = \vect{a}\rfloor A_2 = 0$. Then equation \eqref{eq:ABrrgen} generalizes to 
\begin{align} 
  A\widetilde{B} = A_1 A_2\widetilde{B}
  =& A_1(c_{s+1}c_{s+2}\ldots c_{r-t} \,\,\,+ \,\,\,
  s_{s+1}c_{s+2}\ldots c_{r-t}\bv{i}_{s+1}
        +c_{s+1}s_{s+2}\ldots c_{r-t}\bv{i}_{s+2}+
   \ldots
        +c_{s+1}c_{s+2}\ldots s_{r-t}\bv{i}_{r-t}
  \nonumber \\
  &+
 \ldots+\,\,\,s_{s+1}s_{s+2}\ldots s_{r-t}\bv{i}_{s+1}\bv{i}_{s+2}\ldots\bv{i}_{r-t})
  \,\bv{i}_{r-t+1} \ldots \bv{i}_{r}.
\label{eq:ABrqrgen}
\end{align}
The meaning of parallelity $s$ remains unchanged. Perpendicularity $t$ now means the number of $\pi/2$ principal angles between $\mathsf{A}_2$ and $\mathsf{B}$, which is an invariant, independent of the particular choice of $A_1$. The lowest non-zero grade part of \eqref{eq:ABrqrgen} has dimension $2t+q$ with magnitude  
$\cos\theta_{\mathsf{A}\mathsf{B}}$, while the maximum grade part has magnitude $s_1 \ldots s_{r-t}, \,s_k = \sin \theta_k $. Dividing by the lowest grade part $\cos\theta_{\mathsf{A}\mathsf{B}}\,A_1\bv{i}_{r-t+1} \ldots \bv{i}_{r}$ gives a multivector with maximum grade $2(r-t-s)$, scalar part one, and bivector part
$
t_{s+1}\bv{i}_{s+1}
 +t_{s+2} \bv{i}_{s+2}+\ldots
 +t_{r-t}\bv{i}_{r-t}, 
$
$t_k = \tan \theta_k $. This allows again to compute the principal angles $\theta_k, s<k\leq r-t$.  

\section{Conclusion}

A possible application could be a subspace structure self organizing map
(SOM) type of neural network, mapping the topology of whole data subspace structures faithfully to lower dimensions. The problem of the relative angles between \textit{arbitrary} conformal objects (see \cite{HTBY:Carrier} for definitions) is thus also solved for \textit{all} dimensions. Let $X=S$ be a sphere or $X=F$ be a flat object. If $X\wedge e_{\infty}\neq 0$ (sphere) we set $X=X\wedge e_{\infty}$. Then $X$ always represents the offset embedding flat object. From two such offset embedding flat objects $X,Y$ we can compute the relative angle by setting $A=X\lfloor E$, and $B=Y\lfloor E$, assuming $\mathrm{grade}(X)\geq \mathrm{grade}(Y)$, and applying \eqref{eq:ABrqrgen}. In the future it may be possible to consider the binary (sign) orientation information as well.

\begin{theacknowledgments}
  Soli deo gloria. I do thank my dear family, K. Tachibana, L. Dorst, W. Spr\"{o}ssig and K. G\"{u}rlebeck.
\end{theacknowledgments}

\bibliographystyle{aipproc}

\end{document}